\documentclass[11pt]{article}
\usepackage{amssymb,amsfonts,amsmath, psfrag,eepic,colordvi}
\newtheorem{theo}{Theorem}

\newtheorem{lem} [theo]{Lemma}
\newtheorem{coro}[theo]{Corollary}

\makeatletter \@addtoreset{equation}{section}
\@addtoreset{theo}{}\makeatother

\begin{document}
\begin{center}
{\LARGE\bf
Schr\"oder Paths and Pattern Avoiding Partitions  }\\[16pt]

 Sherry H.F. Yan
 \\[6pt]
Department of Mathematics, Zhejiang Normal University\\
\small Jinhua 321004, P.R. China\\[5pt]
huifangyan@hotmail.com
\end{center}

\noindent {\bf Abstract. } In this paper, we show that
  both $12312$-avoiding partitions  and  $12321$-avoiding  partitions  of the set $[n+1]$  are
 in one-to-one correspondence with  Schr\"oder paths of semilength $n$ without peaks at even level. As
a consequence, the refined enumeration of $12312$-avoiding (resp.
$12321$-avoiding) partitions according to the number of blocks can
be reduced to the enumeration of certain Schr\"oder paths according
to the number of peaks. Furthermore, we get the enumeration of
irreducible $12312$-avoiding (resp. $12321$-avoiding) partitions,
which are closely related to skew Dyck paths.

\vskip 8mm



\vskip 1cm

\section{Introduction and notations}

\begin{figure}[h,t]
\begin{center}
\begin{picture}(200,50)
\setlength{\unitlength}{3mm} \linethickness{0.4pt}

\put(0,0){\circle*{0.2}}
 \put(0,0){\line(1,0){2}}
\put(2,0){\circle*{0.2}}
\put(2,0){\line(1,1){1}}\put(3,1){\circle*{0.2}}\put(3,1){\line(1,1){1}}
\put(4,2){\circle*{0.2}}\put(4,2){\line(1,-1){1}}\put(5,1){\circle*{0.2}}
\put(5,1){\line(1,0){2}}\put(7,1){\circle*{0.2}}\put(7,1){\line(1,0){2}}\put(9,1){\circle*{0.2}}
\put(9,1){\line(1,1){1}}\put(10,2){\circle*{0.2}}\put(10,2){\line(1,1){1}}\put(11,3){\circle*{0.2}}
\put(11,3){\line(1,1){1}}\put(12,4){\circle*{0.2}}
\put(12,4){\line(1,-1){1}}\put(13,3){\circle*{0.2}}
\put(13,3){\line(1,-1){1}}\put(14,2){\circle*{0.2}}
\put(14,2){\line(1,0){2}}\put(16,2){\circle*{0.2}}
\put(16,2){\line(1,-1){1}}\put(17,1){\circle*{0.2}}
\put(17,1){\line(1,0){2}}\put(19,1){\circle*{0.2}}
\put(19,1){\line(1,-1){1}}\put(20,0){\circle*{0.2}}
\put(20,0){\line(1,1){1}}\put(21,1){\circle*{0.2}}
\put(21,1){\line(1,1){1}}\put(22,2){\circle*{0.2}}
\put(22,2){\line(1,-1){1}}\put(23,1){\circle*{0.2}}
\put(23,1){\line(1,-1){1}}\put(24,0){\circle*{0.2}}
\end{picture}
\end{center}
\caption{A UH-free Schr\"{o}der path} \label{fig1}
\end{figure}
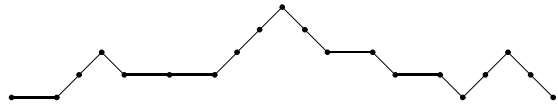

A {\it Schr\"oder path} of semilength $n$ is a lattice path on the
plane from $(0,0)$ to $(2n,0)$ that does not go below the $x$-axis
and consists of up steps $U=(1,1)$, down steps $D=(1,-1)$ and
horizontal steps $H=(2,0)$. They are counted by the larger
Schr\"{o}der numbers (A006318 in \cite{Seq}).  A {\em UH-free}
Schr\"{o}der path is a Schr\"oder path without  up steps   followed
immediately by horizontal steps. A UH-free Schr\"{o}der path of
semilength $12$ is illustrated  as Figure \ref{fig1}.

An up step followed by a down step  in a path is called a {\em
peak}. The {\em level} of an up step (a horizontal step) is defined
as the larger $y$ coordinate of the step. The {\em level} of a peak
is defined as the level of the up step in the peak. Denote by
$\mathcal{SE}_n$ and $\mathcal{SH}_n$ the set of Schr\"{o}der paths
of semilength $n$ without peaks at even level and the set of UH-free
Schr\"{o}der paths of semilength $n$, respectively.

A partition $\pi$ of the set $[n]=\{1,2,\ldots, n\}$ is a collection
$B_1, B_2, \ldots, B_k$ of nonempty disjoint subsets of $[n]$. The
elements of a partition are called blocks. We assume that $B_1, B_2,
\ldots, B_k$ are listed in the increasing order of their minimum
elements, that is $min B_1 <min B_2 <\ldots<min B_k$.  A partition
$\pi$ of $[n]$ with $k$ blocks can also be represented by a sequence
$\pi_1\pi_2\ldots\pi_n$ on the set $\{1, 2, \ldots, k\}$ such that
$\pi_i=j$ if and only if $i\in B_j$. Such a representation is called
the {\em Davenport-Schinzel sequence}  or the {\it canonical
sequential form}. In this paper, we will always represent a
partition by its canonical sequential form.

 In the terminology of
canonical sequential forms, we say that a partition $\pi $ {\em
avoids} a partition $\tau$, or it is $\tau$-{\it avoiding}, if there
is no subsequence which is order-isomorphic to $\tau$ in $\pi$. In
such context, $\tau$ is usually called a pattern. The set of
 $\tau$-avoiding partitions of $[n]$ is denoted $\mathcal{P}_n(\tau)$.
 The enumeration on pattern avoiding
partitions
 has received extensive attention from several authors, see \cite{chen2, Goyt, Jet,  Mansour, Sagan} and  references therein.

By using kernel method,  Mansour and Severini \cite{Mansour} deduced
that the number of $12312$-avoiding partitions of $[n+1]$ is equal
to the number of  Schr\"{o}der paths of semilength $n$ without peaks
at even level (A007317 in \cite{Seq}). Recently, Jelinek and
Mansour\cite{Jet} proved that   the cardinality of
$\mathcal{P}_n(12312)$ is equal to that of $\mathcal{P}_n(12321)$.
In this paper, we will provide  a bijection
 between the set of $12312$-avoiding  partitions of $[n+1]$
and the set of UH-free Schr\"{o}der paths of semilength $n$.  By
making a simple variation of this  bijection, we get a bijection
between the set of $12321$-avoiding  partitions of $[n+1]$ and the
set of UH-free Schr\"{o}der paths of semilength $n$.  A bijection
between the set of UH-free Schr\"{o}der paths of semilength $n$ and
the set of Schr\"{o}der paths of semilength $n$ without peaks at
even level is also provided, which leads to a bijection between
$12312$-avoiding (resp. $12321$-avoiding) partitions of of $[n+1]$
and the set of Schr\"{o}der paths of semilength $n$ without peaks at
even level. As a consequence, the refined enumeration of
 $12312$-avoiding (resp. $12321$-avoiding)
partitions according to the number of blocks can be reduced to the
enumeration of certain Schr\"oder paths according to the number of
peaks. Furthermore, we also get the enumeration of irreducible
$12312$-avoiding (resp. $12321$-avoiding) partitions, which are
closely related to skew Dyck paths.

\section{Bijection between $\mathcal{P}_{n+1}(12312)$ and $\mathcal{SE}_n$}
In this section, we will provide a bijection between the set of
$12312$-avoiding (resp. $12321$-avoiding )partitions of $[n+1]$ and
the set of UH-free Schr\"oder paths of semilength $n$.   A bijection
between the set of UH-free Schr\"oder paths of semilength $n$ and
the set of Schr\"oder paths of semilength $n$ without peaks at even
level is also given, which leads to a bijection between the set of
$12312$-avoiding (resp. $12321$-avoiding )partitions of $[n+1]$ and
the set of Schr\"oder paths of semilength $n$ without peaks at even
level.

Let $\pi$ be a nonempty   partition  of $[n+1]$ with $k$ blocks.
 Then $ \pi $ can be
uniquely decomposed as
\begin{equation}\label{eq.1}
1w_1 2w_2\ldots iw_i\ldots kw_k,
\end{equation}
where $1,2,\ldots, k$ are the {\em left-to-right maxima} of $\pi$
and each $w_i$ is a possibly empty word on $[i]$.
 For, $1\leq i\leq k$, denote by $w_i\setminus \{i\}$
the word  obtained from $w_i$ by deleting all the $i's$. The
following property of 12312-avoiding (resp. 12321-avoiding )
partitions can be verified easily and we omit the proof here.
\begin{lem}\label{lem1}
A partition $\pi$ is 12312-avoiding (resp. 12321-avoiding) partition
with $k$ blocks if and only if  the word $w_1\setminus
\{1\}w_2\setminus \{2\}\ldots w_k\setminus \{k\}$ is in weakly
decreasing (increasing) order.
\end{lem}

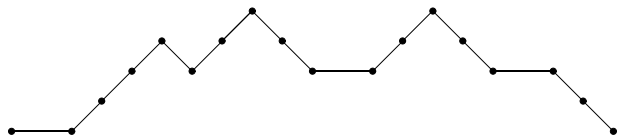
\begin{figure}[h,t]
\begin{center}
\begin{picture}(200,70)
\setlength{\unitlength}{4mm} \linethickness{0.4pt}

\put(0,0){\circle*{0.2}}
 \put(0,0){\line(1,0){2}}
 \put(2,0){\circle*{0.2}}
 \put(2,0){\line(1,1){1}}
\put(3,1){\circle*{0.2}}
 \put(3,1){\line(1,1){1}}
 \put(4,2){\circle*{0.2}}
 \put(4,2){\line(1,1){1}}
 \put(5,3){\circle*{0.2}}
 \put(5,3){\line(1,-1){1}}
 \put(6,2){\circle*{0.2}}
 \put(6,2){\line(1,1){1}}
 \put(7,3){\circle*{0.2}}
 \put(7,3){\line(1,1){1}}
 \put(8,4){\circle*{0.2}}
 \put(8,4){\line(1,-1){1}}
  \put(9,3){\circle*{0.2}}
 \put(9,3){\line(1,-1){1}}
  \put(10,2){\circle*{0.2}}
 \put(10,2){\line(1,0){2}}
  \put(12,2){\circle*{0.2}}
 \put(12,2){\line(1,1){1}}
  \put(13,3){\circle*{0.2}}
 \put(13,3){\line(1,1){1}}
  \put(14,4){\circle*{0.2}}
 \put(14,4){\line(1,-1){1}}
  \put(15,3){\circle*{0.2}}
 \put(15,3){\line(1,-1){1}}
  \put(16,2){\circle*{0.2}}
 \put(16,2){\line(1,0){2}}
  \put(18,2){\circle*{0.2}}
 \put(18,2){\line(1,-1){1}}
  \put(19,1){\circle*{0.2}}
 \put(19,1){\line(1,-1){1}}
  \put(20,0){\circle*{0.2}}
 \put(7,3){\line(1,1){1}}
\end{picture}
\end{center}
\caption{A UH-free Schr\"oder path  of semilength $10$.}
\label{UH-free}
\end{figure}

 Now, we proceed to construct a map $\sigma $ from
$\mathcal{P}_{n+1}(12312)$ to $\mathcal{SH}_n$. Given a
$12312$-avoiding partition $\pi$ of $[n+1]$ with $k$ blocks, if
$\pi=1$, then let $\sigma(\pi)$ be the empty path. Otherwise,
suppose that $\pi$ is decomposed as (\ref{eq.1}) and for
$i=1,2,\ldots, k-1$, denote by $d_i$ the number of occurrences of
$i$ which are right to the first occurrence of $i+1$. We read the
decomposition from left to right and generate a path $\sigma(\pi)$
as follows: when a left-to-right maximum $i$ ($i\geq 2$) is read, we
adjoin  $d_{i-1}+1$ successive up steps followed by one down step;
when each element less than $i$ in any word $w_i$ ($1\leq i\leq k$)
is read, we adjoin one down step; when  each element $i$ in any word
$w_i$ ($1\leq i\leq k$) is read, we adjoin one horizontal step.
Lemma \ref{lem1} ensures that the obtained path $\sigma(\pi)$ is a
well defined UH-free Schr\"oder path  of semilength $n$. For
instance, a $12312$-avoiding partition $\pi=11232343411$ of $[11]$
can be decomposed as $1w_12w_23w_34w_4$, where $w_1=1, w_3=23,
w_4=3411, d_1=2, d_2=1 , d_3=1,$ and $w_2$ is empty.
 The corresponding  UH-free Schr\"oder path  $\sigma(\pi)$ of
semilength $10$  is illustrated as Figure \ref{UH-free}.

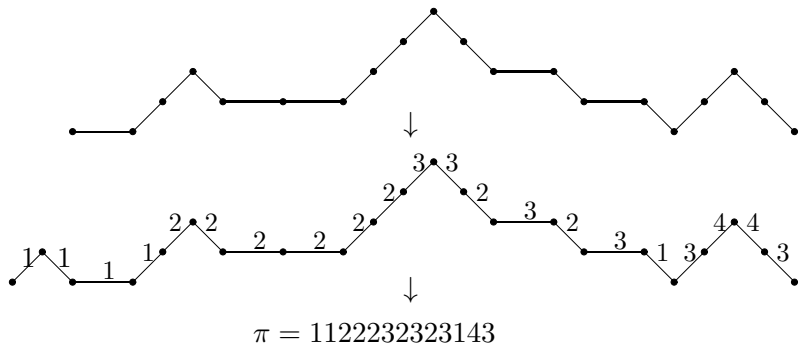
\begin{figure}[h,t]
\begin{center}
\begin{picture}(300,120)
\setlength{\unitlength}{4mm} \linethickness{0.4pt}

\put(2,1){\circle*{0.2}}
 \put(2,1){\line(1,0){2}}
 \put(4,1){\circle*{0.2}}
 \put(3,1.1){\small$1$}
\put(4,1){\line(1,1){1}}\put(5,2){\circle*{0.2}}
\put(4.3,1.7){\small$1$}
 \put(5,2){\line(1,1){1}}
\put(6,3){\circle*{0.2}} \put(5.2,2.7){\small$2$}

\put(6,3){\line(1,-1){1}}\put(7,2){\circle*{0.2}}
\put(6.4,2.7){\small$2$}
\put(7,2){\line(1,0){2}}\put(9,2){\circle*{0.2}}
 \put(8,2.1){\small$2$}
\put(9,2){\line(1,0){2}}\put(11,2){\circle*{0.2}}
 \put(10,2.1){\small$2$}
\put(11,2){\line(1,1){1}}\put(12,3){\circle*{0.2}}
\put(11.3,2.7){\small$2$}
\put(12,3){\line(1,1){1}}\put(13,4){\circle*{0.2}}
\put(12.3,3.7){\small$2$}
\put(13,4){\line(1,1){1}}\put(14,5){\circle*{0.2}}
\put(13.3,4.7){\small$3$}
\put(14,5){\line(1,-1){1}}\put(15,4){\circle*{0.2}}
\put(14.4,4.7){\small$3$}
\put(15,4){\line(1,-1){1}}\put(16,3){\circle*{0.2}}
\put(15.4,3.7){\small$2$}
\put(16,3){\line(1,0){2}}\put(18,3){\circle*{0.2}}
\put(17,3.1){\small$3$}
\put(18,3){\line(1,-1){1}}\put(19,2){\circle*{0.2}}
\put(18.4,2.7){\small$2$}
\put(19,2){\line(1,0){2}}\put(21,2){\circle*{0.2}}
\put(20,2.1){\small$3$}
\put(21,2){\line(1,-1){1}}\put(22,1){\circle*{0.2}}
\put(21.4,1.7){\small$1$}
\put(22,1){\line(1,1){1}}\put(23,2){\circle*{0.2}}
\put(22.3,1.7){\small$3$}
\put(23,2){\line(1,1){1}}\put(24,3){\circle*{0.2}}
\put(23.3,2.7){\small$4$}
\put(24,3){\line(1,-1){1}}\put(25,2){\circle*{0.2}}
\put(24.4,2.7){\small$4$}
\put(25,2){\line(1,-1){1}}\put(26,1){\circle*{0.2}}
\put(25.4,1.7){\small$3$} \put(0,1){\circle*{0.2}}
\put(0,1){\line(1,1){1}}\put(1,2){\circle*{0.2}}
\put(1,2){\line(1,-1){1}} \put(0.3, 1.5){$1$} \put(1.5, 1.5){$1$}
\put(2,6){\circle*{0.2}}
 \put(2,6){\line(1,0){2}}
 \put(4,6){\circle*{0.2}}
\put(4,6){\line(1,1){1}}\put(5,7){\circle*{0.2}}
 \put(5,7){\line(1,1){1}}
\put(6,8){\circle*{0.2}}
\put(6,8){\line(1,-1){1}}\put(7,7){\circle*{0.2}}
\put(7,7){\line(1,0){2}}\put(9,7){\circle*{0.2}}
\put(9,7){\line(1,0){2}}\put(11,7){\circle*{0.2}}
 \put(11,7){\line(1,1){1}}\put(12,8){\circle*{0.2}}
\put(12,8){\line(1,1){1}}\put(13,9){\circle*{0.2}}
\put(13,9){\line(1,1){1}}\put(14,10){\circle*{0.2}}
\put(14,10){\line(1,-1){1}}\put(15,9){\circle*{0.2}}
\put(15,9){\line(1,-1){1}}\put(16,8){\circle*{0.2}}
\put(16,8){\line(1,0){2}}\put(18,8){\circle*{0.2}}
\put(18,8){\line(1,-1){1}}\put(19,7){\circle*{0.2}}
\put(19,7){\line(1,0){2}}\put(21,7){\circle*{0.2}}
\put(21,7){\line(1,-1){1}}\put(22,6){\circle*{0.2}}
\put(22,6){\line(1,1){1}}\put(23,7){\circle*{0.2}}
\put(23,7){\line(1,1){1}}\put(24,8){\circle*{0.2}}
\put(24,8){\line(1,-1){1}}\put(25,7){\circle*{0.2}}
\put(25,7){\line(1,-1){1}}\put(26,6){\circle*{0.2}}
\put(13, 6){$\downarrow$}
\put(13, 0.5){$\downarrow$}
\put(8,-1){$\pi=1122232323143$}
\end{picture}
\end{center}
\caption{ An example of the reverse map of $\sigma$.} \label{12312}
\end{figure}

\begin{figure}[h,t]
\begin{center}
\begin{picture}(300,120)
\setlength{\unitlength}{4mm} \linethickness{0.4pt}

\put(2,1){\circle*{0.2}}
 \put(2,1){\line(1,0){2}}
 \put(4,1){\circle*{0.2}}
 \put(3,1.1){\small$1$}
\put(4,1){\line(1,1){1}}\put(5,2){\circle*{0.2}}
\put(4.3,1.7){\small$1$}
 \put(5,2){\line(1,1){1}}
\put(6,3){\circle*{0.2}} \put(5.2,2.7){\small$2$}

\put(6,3){\line(1,-1){1}}\put(7,2){\circle*{0.2}}
\put(6.4,2.7){\small$2$}
\put(7,2){\line(1,0){2}}\put(9,2){\circle*{0.2}}
 \put(8,2.1){\small$2$}
\put(9,2){\line(1,0){2}}\put(11,2){\circle*{0.2}}
 \put(10,2.1){\small$2$}
\put(11,2){\line(1,1){1}}\put(12,3){\circle*{0.2}}
\put(11.3,2.7){\small$2$}
\put(12,3){\line(1,1){1}}\put(13,4){\circle*{0.2}}
\put(12.3,3.7){\small$2$}
\put(13,4){\line(1,1){1}}\put(14,5){\circle*{0.2}}
\put(13.3,4.7){\small$3$}
\put(14,5){\line(1,-1){1}}\put(15,4){\circle*{0.2}}
\put(14.4,4.7){\small$3$}
\put(15,4){\line(1,-1){1}}\put(16,3){\circle*{0.2}}
\put(15.4,3.7){\small$1$}
\put(16,3){\line(1,0){2}}\put(18,3){\circle*{0.2}}
\put(17,3.1){\small$3$}
\put(18,3){\line(1,-1){1}}\put(19,2){\circle*{0.2}}
\put(18.4,2.7){\small$2$}
\put(19,2){\line(1,0){2}}\put(21,2){\circle*{0.2}}
\put(20,2.1){\small$3$}
\put(21,2){\line(1,-1){1}}\put(22,1){\circle*{0.2}}
\put(21.4,1.7){\small$2$}
\put(22,1){\line(1,1){1}}\put(23,2){\circle*{0.2}}
\put(22.3,1.7){\small$3$}
\put(23,2){\line(1,1){1}}\put(24,3){\circle*{0.2}}
\put(23.3,2.7){\small$4$}
\put(24,3){\line(1,-1){1}}\put(25,2){\circle*{0.2}}
\put(24.4,2.7){\small$4$}
\put(25,2){\line(1,-1){1}}\put(26,1){\circle*{0.2}}
\put(25.4,1.7){\small$3$} \put(0,1){\circle*{0.2}}
\put(0,1){\line(1,1){1}}\put(1,2){\circle*{0.2}}
\put(1,2){\line(1,-1){1}} \put(0.3, 1.5){$1$} \put(1.5, 1.5){$1$}
\put(2,6){\circle*{0.2}}
 \put(2,6){\line(1,0){2}}
 \put(4,6){\circle*{0.2}}
\put(4,6){\line(1,1){1}}\put(5,7){\circle*{0.2}}
 \put(5,7){\line(1,1){1}}
\put(6,8){\circle*{0.2}}
\put(6,8){\line(1,-1){1}}\put(7,7){\circle*{0.2}}
\put(7,7){\line(1,0){2}}\put(9,7){\circle*{0.2}}
\put(9,7){\line(1,0){2}}\put(11,7){\circle*{0.2}}
 \put(11,7){\line(1,1){1}}\put(12,8){\circle*{0.2}}
\put(12,8){\line(1,1){1}}\put(13,9){\circle*{0.2}}
\put(13,9){\line(1,1){1}}\put(14,10){\circle*{0.2}}
\put(14,10){\line(1,-1){1}}\put(15,9){\circle*{0.2}}
\put(15,9){\line(1,-1){1}}\put(16,8){\circle*{0.2}}
\put(16,8){\line(1,0){2}}\put(18,8){\circle*{0.2}}
\put(18,8){\line(1,-1){1}}\put(19,7){\circle*{0.2}}
\put(19,7){\line(1,0){2}}\put(21,7){\circle*{0.2}}
\put(21,7){\line(1,-1){1}}\put(22,6){\circle*{0.2}}
\put(22,6){\line(1,1){1}}\put(23,7){\circle*{0.2}}
\put(23,7){\line(1,1){1}}\put(24,8){\circle*{0.2}}
\put(24,8){\line(1,-1){1}}\put(25,7){\circle*{0.2}}
\put(25,7){\line(1,-1){1}}\put(26,6){\circle*{0.2}}\put(13,
6){$\downarrow$}
\put(13, 0.5){$\downarrow$} \put(8,-1){$\pi=1122231323243$}

\end{picture}
\end{center}
\caption{An example of the reverse map of $\phi$.} \label{12321}
\end{figure}
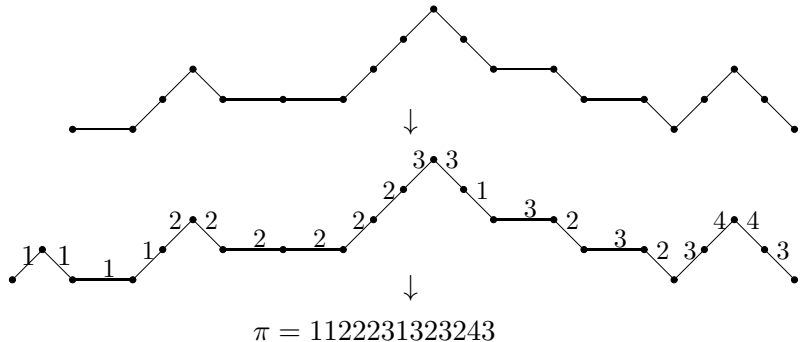

Conversely, we can get a $12312$-avoiding partition of $[n+1]$  from
a UH-free Schr\"oder path $P$ of semilength $n$. If $P$ is empty,
then let $\sigma^{-1}(P)=1$, otherwise suppose that $P$ has $k$
peaks. Then we can get a word $\sigma^{-1}(P)$ as the following
procedure.
\begin{itemize}
 \item[Step 1.] Firstly, add a peak at the very begining of $P$ and denote by $P'$ the obtained path;
 \item[Step 2.] Secondly, label all the  up steps in peaks of $P'$ with the
alphabet $\{1,2,\ldots, k+1\}$ from left to right and label each
remaining up step $s$ and  each  horizontal step $h$ with the
maximum alphabet which are left to the steps $s$ and $h$,
respectively;

\item [Step 3.] Thirdly,  if  a down step $s$ is  in a peak, label $s$ with the
same label as the label of the  up step in the same peak; Otherwise,
suppose that $L^{U}$ (resp. $L^{D}$) is the multiset of all the
labels of the up (resp. down) steps left to the   step $s$. Then
label $s$ with the maximum element of the multiset obtained from
$L^{U}$ by removing all the elements of $L^{D}$;
\item[Step 4.] Lastly, let $\sigma^{-1}(P)$ be a word obtained by reading
the labels of all the down steps and horizontal steps of $P'$
successively.

\end{itemize}
Obviously, the obtained word $\sigma^{-1}(P)$ is a $12312$-avoiding
partition of $[n+1]$. An example of the reverse map of $\sigma$ is
shown in Figure \ref{12312}.
\begin{theo}\label{sigma}
The map $\sigma$ is a bijection between the set of $12312$-avoiding
partitions of $[n+1]$ and the set of UH-free Schr\"oder paths of
semilength $n$.
\end{theo}
 We define a map $\phi$ from
$\mathcal{P}_{n+1}(12321)$ to $\mathcal{SH}_n$ the same as the map
$\sigma$ and define the reverse of $\phi$ the same as the reverse of
$\sigma$ except that in Step 3 we label the down step $s$ not in a
peak by the minimum element of the multiset obtained from $L^{U}$ by
removing all the elements of $L^{D}$. It is easy to check that
$\phi$ is a bijection between the set of $12321$-avoiding partitions
of $[n+1]$ and the set of UH-free Schr\"oder paths of semilength
$n$. An example of the reverse map of $\phi$ is illustrated as
Figure \ref{12321}.

\begin{theo}\label{phi}
The map  $\phi$ is a bijection between  UH-free Schr\"oder paths of
semilength $n$ and $\mathcal{P}_{n+1}(12321)$.
\end{theo}

In order to get a bijection between $\mathcal{P}_{n+1}(12312)$ and
$\mathcal{SE}_n$, we should provide a bijection between
$\mathcal{SH}_n$ and $\mathcal{SE}_n$. Now we proceed to construct
the map $\psi$ from $\mathcal{SH}_n$ and $\mathcal{SE}_n$. Given a
UH-free Schr\"oder path $P\in \mathcal{SH}_n$, if it is empty, then
let $\psi(P)$ be an empty path. Otherwise, we can  get $\psi(P)$
recursively as follows:
\begin{itemize}
\item If $P=HP'$, then let $\psi(P)=H\psi(P')$, where $P'$ is a
possibly empty UH-free Schr\"oder path;
\item If $P=UDP'$, then let $\psi(P)=UD\psi(P')$, where $P'$ is a
possibly empty UH-free Schr\"oder path;
\item If $P=U^kDP_1DP_2\ldots DP_k$, where $k\geq 2$, $U^k$ denotes  $k$ consecutive
up steps and  for $1\leq i\leq k$, each $P_i$  is a possibly empty
UH-free Schr\"oder path, then let $\psi(P)=UP'_1P'_2\ldots
P'_{k-1}D\psi(P_k)$ such that for $1\leq i\leq k-1$, each $
P'_{i}=H$ if $P_i$ is empty and $P'_i=U\psi(P_i)D$, otherwise.
\end{itemize}
Obviously, the obtained path $\psi(P)$ is a Schr\"oder path of
semilength $n$ without peaks at even level. It is easy to check that
the map $\psi$ is reversible. For the convenience of simplicity, we
omit the reverse map of $\psi$.
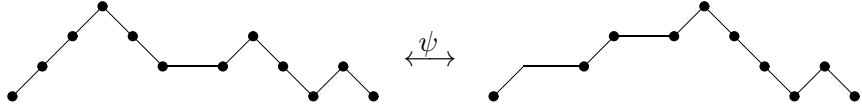
\begin{figure}[h,t]
\begin{center}
\begin{picture}(300,70)
\setlength{\unitlength}{4mm} \linethickness{0.4pt}
\put(0,0){\circle*{0.3}}\put(0,0){\line(1,1){1}}
\put(1,1){\circle*{0.3}}\put(1,1){\line(1,1){1}}
\put(2,2){\circle*{0.3}}\put(2,2){\line(1,1){1}}
\put(3,3){\circle*{0.3}}\put(3,3){\line(1,-1){1}}
\put(4,2){\circle*{0.3}}\put(4,2){\line(1,-1){1}}
\put(5,1){\circle*{0.3}}\put(5,1){\line(1,0){2}}
\put(7,1){\circle*{0.3}}\put(7,1){\line(1,1){1}}
\put(8,2){\circle*{0.3}}\put(8,2){\line(1,-1){1}}
\put(9,1){\circle*{0.3}}\put(9,1){\line(1,-1){1}}
\put(10,0){\circle*{0.3}}\put(10,0){\line(1,1){1}}
\put(11,1){\circle*{0.3}}\put(11,1){\line(1,-1){1}}
\put(12,0){\circle*{0.3}}
\put(13,1){$\longleftrightarrow$} \put(13.5,1.5){$\psi$}
\put(16,0){\circle*{0.3}}\put(16,0){\line(1,1){1}}
\put(17,1){\line(1,0){2}}\put(19,1){\circle*{0.3}}
\put(19,1){\line(1,1){1}}\put(20,2){\circle*{0.3}}\put(20,2){\line(1,0){2}}
\put(22,2){\circle*{0.3}}\put(22,2){\line(1,1){1}}
\put(23,3){\circle*{0.3}}\put(23,3){\line(1,-1){1}}
\put(24,2){\circle*{0.3}}\put(24,2){\line(1,-1){1}}
\put(25,1){\circle*{0.3}}\put(25,1){\line(1,-1){1}}
\put(26,0){\circle*{0.3}}\put(26,0){\line(1,1){1}}
\put(27,1){\circle*{0.3}}\put(27,1){\line(1,-1){1}}
\put(28,0){\circle*{0.3}}
\end{picture}
\end{center}
\caption{A UH-free Schr\"oder path and its corresponding Schr\"oder
path without peaks at even level.} \label{uh}
\end{figure}

\begin{theo}\label{sch}
The map $\psi$ is a bijection between the set of UH-free Schr\"oder
paths of semilength $n$ and the set of Schr\"oder path of semilength
$n$ without peaks at even level.
\end{theo}
Combining Theorems \ref{sigma},  \ref{phi} and \ref{sch}, we have
the following results.
\begin{theo}
The map $\psi\cdot \sigma $ (resp. $\psi\cdot \phi $)  is a
bijection between the set of $12312$-avoiding (resp.
$12321$-avoiding) partitions of $[n+1]$ and the set of Schr\"oder
paths of semilength $n$ without peaks at even level.
\end{theo}

\section{Refined enumerations}
In this section, we aim to  get the refined enumeration of
$12312$-avoiding (resp. $12321$-avoidng )partitions according to the
number of blocks.
  By restricting the peaks in
Schr\"oder paths, we get the enumeration of irreducible
$12321$-avoiding (resp. 12321-avoiding)  partitions. From the
construction of the maps $\sigma$ and $\phi$, we get that each block
apart from the first block in a 12312-avoiding (resp.
12321-avoiding) partition $\pi$ brings up a peak in its
corresponding UH-free Schr\"oder path $\sigma(\pi)$ (resp.
$\phi(\pi)$). Hence, we get the following result.
\begin{coro}\label{coro1}
Let $\pi$ be  a 12312-avoiding (resp. 12321-avoiding) partition on
$[n+1]$ with $k+1$ blocks, then  $\sigma(\pi)$ (resp. $\phi(\pi)$)
is  a UH-free Schr\"oder path of semilength $n$ with $k$ peaks.
\end{coro}
A {\em Dyck path} of semilength $n$ is a lattice path on the plane
from $(0,0)$ to $(2n,0)$ that does not go below the $x$-axis and
consists of up steps $U=(1,1)$ and down steps $D=(1,-1)$. The number
of Dyck paths of semilength $n$  with $k$ peaks is counted by the
Narayana number $$N_{n,k}={1\over n}{n\choose k}{n\choose k-1}.$$
Note that any UH-free Schr\"oder path of semilength $n$  with $k$
peaks can be obtained from a Dyck path of semilength $j$ $(0\leq
j\leq n)$  with $k$ peaks by inserting $n-j$ horizontal steps into
the positions after down steps and the position at the very
beginning of the Dyck path. The number of such arrangement is equal
to ${n\choose j}$. Hence, the number of UH-free Schr\"oder paths of
semilength $n$  with $k$  peaks is counted by $\sum_{j=k}^{n}{1\over
j}{j\choose k}{j\choose k-1}{n\choose j}$. From Corollary
\ref{coro1}, we get the following result.

 \begin{coro}
 The number of 12312-avoiding (resp. 12321-avoiding ) partitions on $[n+1]$ with $k+1$ ($k\geq 1$) blocks is equal to
 $$
 \sum_{j=k}^{n}{1\over j}{j\choose k-1}{j\choose k} {n\choose j}.
 $$

 \end{coro}

 A partition $P$ of $[n]$ is called an {\em irreducible
partition} if for any $m\in [n-1]$,    $P$ can not be reduced to two
smaller partitions $P_1$ and $P_2$ such that $P_1$ is  a partition
of $[m]$ and $P_2$ is a partition of $\{m+1, m+2, \ldots, n\}$.
Irreducible partitions have been
 studied by Lehner \cite{Le}.  In fact, a
 partition $\pi$ of $[n]$ is irreducible if and only if for any
 element $i\in [n]$, there is at least one occurrence of an element
 $j$ which is less than $i$ and right to the first occurrences of
 $i$. Hence, by the construction of the maps $\sigma$ and $\phi$, we see that if $\pi$ is
 irreducible, then its corresponding Schr\"{o}der
path $\sigma(\pi)$ (resp. $\phi(\pi)$) has no peaks at level one.
\begin{coro}
The map $\sigma$ (resp. $\phi$) is a bijection between the set
 of irreducible
$12312$-avoiding (resp. $12321$-avoiding)  partitions on $[n+1]$ and
the set of UH-free Schr\"{o}der paths of semilength $n$ without
peaks at level one.
\end{coro}

Denote by $\mathcal{SH'}_n$ the set of   UH-free Schr\"{o}der paths
of semilength $n$ without peaks at level one. Let $s_n$ and $s'_n$
the cardinality of $\mathcal{SH}_n$ and  $\mathcal{SH'}_n$,
respectively. Let $f(x)=\sum_{n=0}^{\infty}s_nx^n$  and
$f'(x)=\sum_{n=0}^{\infty}s'_nx^n$ where $s_0=1 $ and $s'_0=1$.
Then, it is easy to get the following recurrence relations:
$$
f(x)=1+2xf(x)+xf(x)(f(x)-1-xf(x)),
$$
and
$$
f'(x)=1+xf'(x)+xf'(x)(f(x)-1-xf(x)).
$$
Hence, we have
$$
f(x)={1-x-\sqrt{1-6x+5x^2}\over 2(x-x^2)},
$$
and
$$
f'(x)={1\over 1-x(1-x)f(x)}={2\over 1+x+\sqrt{1-6x+5x^2} },
$$
which is the generating function for skew Dyck paths of length $n$
ending with a down step, see \cite[A033321]{Seq}. A {\em skew Dyck
path} is a path in the first quadrant which begins at the origin,
ends on the $x$-axis, consists of steps $U=(1,1)$, $D=(1,-1)$,  and
$L=(-1,-1)$ so that never lie below the x-axis and  up and left
steps do not overlap.
\begin{coro}
The number of irreducible  $12312$-avoiding (resp. $12321$-avoiding)
partitions of $[n+1]$ is equal to the number of skew Dyck paths of
semilength $n$ ending with a down step.
\end{coro}

 \vskip 5mm

\small

\end{document}